\documentclass[12pt]{article}
\usepackage{amsmath}
\usepackage{amssymb}
\usepackage{mathptmx}
\usepackage{amsfonts}
\usepackage[mathscr]{eucal}
\usepackage[dvips]{graphicx}
\usepackage{color}

\newcommand{\Yukawa}{\textrm{Y}}
\newcommand{\Dirac}{\textrm{D}}
\newcommand{\KG}{ \textrm{KG} }
\newcommand{\I}{\textrm{I}}
\newcommand{\bos}{\textrm{b}}
\newcommand{\fer}{\textrm{f}}
\newcommand{\fin}{\textrm{fin}}
\newcommand{\Dir}{ \textrm{D} }
\newcommand{\ess}{ \textrm{ess} }

\newcommand{\stateKG}{L^{2}( \mathbf{R}^{3}_{\mathbf{k}}  )}
\newcommand{\stateDirac}{L^{2}( \mathbf{R}^{3}_{\mathbf{p}}  ; \mathbf{C}^{4})}

\newcommand{\Fb}{\mathscr{F}_{\textrm{b}}  } 
\newcommand{\Ff}{\mathscr{F}_{\textrm{f}}  }

\newcommand{\HDirac}{H_{\textrm{D}}}
\newcommand{\HKG}{H_{\textrm{KG}}}
\newcommand{\HYukawa}{H_{\textrm{Y}}}

\newcommand{\HI}{H_{\textrm{I}}}

\newcommand{\Hzero}{H_{0}}

\newcommand{\mbf}[1]{\ensuremath{\mathbf{#1}}}
\newcommand{\ms}[1]{\ensuremath{\mathscr{#1}}}

\newcommand{\sqzb}[1]{\ensuremath{d\Gamma_{\textrm{b}}({#1}) }}
\newcommand{\sqzf}[1]{\ensuremath{d\Gamma_{\textrm{f}}({#1}) }}

\newcommand{\tens}{\otimes}

\newcommand{\nstens}{\otimes^{n}_{s}}
\newcommand{\natens}{\otimes^{n}_{a}}

\newcommand{\eltwo}{L^{2}(\mathbf{R}^{3})}
\newcommand{\Rthree}{\mathbf{R}^{3} }
\newcommand{\dx}{d \mathbf{x} }

\newcommand{\dk}{d \mathbf{k} }

\newcommand{\restr}{\upharpoonright}

\newcommand{\Rd}{\mathbf{R}^{d} }

\newcommand{\cutoff}[2]{\ensuremath{
\chi_{\textrm{#1}} (\mathbf{#2} )}}

\newtheorem{theorem}{Theorem}[section]

\newtheorem{lemma}[theorem]{Lemma}
\newtheorem{corollary}[theorem]{Corollary}

\begin{document}
\begin{center}
{\LARGE  Essential Spectrum of  a Fermionic \\ Quantum 
Feild Model} \\
 $\;$ \\
 Toshimitsu Takaesu \\
 $\;$ \\
\textit{Faculty of Science and Engineering, Gunma University,\\ Gunma, 376-8515, Japan }
\end{center}

\begin{quote}
\textbf{Abstract}. An interaction   system of  a fermionic quantum field is considered. The state space is defined  by a tensor product space of a fermion Fock space and a Hilbert space. It is  assumed that  the total Hamiltonian is a self-adjoint operator on the state space and bounded from below. Then it is proven that  a subset 
of real numbers is  the essential spectrum of the total Hamiltonian. It is  applied to the  system of a Dirac field coupled to a Klein-Gordon field.  Then the HVZ theorem for the system  is obtained. \\
$\; $ \\
{\small
key words : Essential Spectrum, Fock Space, Quantum Field Theory}.
\end{quote}
 
\section{Introduction and an Main Theorem }
In this paper, an interaction system of a fermionic field is invested. Let $\ms{K}$ be a Hilbert space and $\Ff ( \ms{K} )$ the fermion Fock space over  $\ms{K}$. Let  $\ms{T}$ be a 
Hilbert pace. The state space  of the interaction system is defined by  
\begin{equation}
\ms{H} \; = \; \Ff ( \ms{K} ) \tens \ms{T} .
\end{equation}
Let $K$ be a self-adjoint operator on $\ms{K}$ with ker $K = \{ 0\}$ and $\sqzf{K} $  the second quantization of $K$. 
Let $T$ be a self-adjoint  operator on $\ms{T}$. We assume that $K$ is non-negative and $T$ is bounded from below. 
The free Hamiltonian is given by 
 \begin{equation}
 \Hzero \; = \; \sqzf{K } \tens I \; + \; I \tens  T .
 \end{equation}
It is seen that   $H_{0} $ is  self-adjoint  on 
$\ms{D}(\Hzero ) = \ms{D} (\sqzf{K } \tens I ) \cap \ms{D} (I \tens T)$ and bounded from below. The total Hamiltonian is given by
 \begin{equation}
H \; = \; H_0 \; +  \HI , 
\end{equation}
where $\HI$ is a symmetric operator on $\ms{H}$.

$\; $ \\ 
We are interested in   the essential spectrum of $H$. 
Locations of  essential spectrum of quantum filed Hamiltonians  are  investigated by  methods in scattering theory  in \cite{Am04,DG99} and by a weak commutator method in  \cite{Arai00}.
 In this paper, we apply the weak commutator method in \cite{Arai00}, which is mentioned below.   
Let $X$ and $Y$ be  densely defined linear operators on a Hilbert space $\ms{X}$. Assume that there exist a linear operator  $Z$ and a  dense subspace $\ms{M}$ such that $\ms{M} \subset \ms{D} (Z)  \cap \ms{D} (X)  \cap \ms{D} (X^{\ast})  \cap  \ms{D} (Y) \cap \ms{D} (Y^{\ast})  $ and for  all $\Phi, \Psi  \in \ms{M}$, 

\[
\qquad \quad 
 ( X^\ast \Phi , Y \Psi  )  - (Y^{\ast} \Phi , X \Psi ) \; = \;  (\Phi , Z \Psi ) 
 .
\]
Then the  restriction of $Z$  to $\ms{M}$
  is called a   weak commutator  of $X$ and $Y$ on $\ms{M}$, and  denoted by  
$     [X , Y ]^0_{\ms{M}}  $.

$\;$ \\
We  suppose   conditions below : 
\begin{quote}
\textbf{(A.1)} $H$ is self-adjoint on $\ms{D} (H) = \ms{D} (H_0) \cap \ms{D} (\HI )$ and bounded from below. \\
$\; $ \\
\textbf{(A.2)} For all $h \in \ms{D} (K)$, $[\HI  , B (h) \tens I ]^0_{\ms{D}(H)} $ and  $[\HI  , B^\ast (h) \tens I ]^0_{\ms{D}(H)} $ exist, where  $B(h) $ and $B^{\ast} (h)$ denote the annihilation operator and the creation operator, respectively. In addition, for   any sequence
 $ \{ h_n \}_{n=1}^\infty $  of $\ms{D} (K) $ such that w-$\lim\limits_{n \to \infty} h_n = 0$ and  $\| h_n \| = 1 $, $n \geq 1$,  it follows that for all $\Psi \in \ms{D} (H)$, 
\begin{align*}
&\textrm{(1)} \quad \textrm{s-}\lim_{n \to \infty}[\HI  , B (h_n) \tens I ]^0_{\ms{D}(H)} \, \Psi   = 0 , \qquad \\
&\textrm{(2)} \quad \textrm{s-}\lim_{n \to \infty}[\HI  , B^\ast (h_n)  \tens I ]^0_{\ms{D}(H)} \, \Psi   = 0 . \qquad 
\end{align*}
\end{quote}

$\;$ \\
For a self-adjoint operator $X$, the spectrum of $X$ is denoted by   $\sigma (X)$, and the essential spectrum by $\sigma_{\ess} (X)$
From \textbf{(A.1)},  it follows that  
 $ E_{0} (H) > - \infty $, where $E_{0} (H) = \inf \sigma (H) $ is the ground state energy of $H$. 

\begin{theorem}  \label{main}
Assume \textbf{(A.1)} and \textbf{(A.2)}. Then 
\[
\overline{  \{  E_0 ( H )  + \lambda \left| \frac{}{} \right.  \lambda \in  \sigma_{\ess} (K) \backslash \{ 0  \} \}}      \; \subset \sigma_{\ess} (H) , 
\]
where $\overline{J} $ denotes the closure of $J \subset \mbf{R}$.
\end{theorem}
$\;$ \\
In section 3 we consider an application of  Theorem \ref{main} to a system in the Yukawa theory. The Yukawa theory  describes  systems of  fermionic  fields coupled to  bosonic fields (e.g., \cite{DP12, Fr73, GJ85}). We consider the system of a Dirac field coupled a Klein-Gordon field. The total Hamiltonian is a defined on a boson-fermion Fock space. The existence of a positive spectral gap above the ground sate energy is proven in \cite{Ta11}. By  Theorem \ref{main} and results in \cite{Arai00,Ta11}, it is proven that a subset of $\mbf{R}$ is  equal to  the essential spectrum of the total Hamiltonian. This result can be regarded as a  quantum field version of the HVZ theorem for quantum mechanics systems \cite{CFKS}.

$\;$ \\
This paper is organized as follows. 
In Section 2, basic properties of Fock spaces are explained  and  the proof of  Theorem \ref{main} is given. In section 3, an application of Theorem \ref{main} to the Yukawa model is considered.

\section{Proof of Theorem \ref{main}}
\subsection{Fermion Fock Space and Boson Fock Space}
First we introduce Fermion Fock space (e.g., \cite{Gu12,Tha}).
The fermion Fock space over a Hilbert space $\ms{X}$ is defined by 
 $ \Ff (\ms{X}) = \oplus_{n=0}^{\infty}
 (  \natens \ms{X} ) \, $, where $\natens \ms{X} $ denotes the $n$-fold anti-symmetric tensor product of $\ms{X} $ with
$\; \tens_{a}^{0} \ms{X} := \mbf{C}$.
The Fock vacuum  is defined by
$ \Omega_{\fer}  = \{   1, 0, 0, \cdots \} \, \in \, \Ff (\ms{X}) $. 
The annihilation  operator is denoted by $B(f)$, $f \in \ms{X}$ and  
the  creation operator by $\; B^{\ast}(g)$, $ g \in \ms{X}$.  
For a dense subspace $\ms{M} \subset \ms{X} $, the finite particle subspace $\Ff^{\fin} (\ms{M})$ is  the linear hull of vectors of the form $ \Psi =  B^{\ast } (f_1) \cdots  B^{\ast } (f_n) \Omega_{\fer}$, $ f_j \in \ms{M}$, $j=1, \cdots n$, $n \in \mbf{N} $.   
It is known that
$B(f)$ and $B^{\ast } (g) $ is bounded  with 
\begin{equation}
\| B(f) \| = \| f  \| , \quad \| B^{\ast } (g) \| = \| g \| , \label{boundBBast}
\end{equation}
 respectively. They satisfy   
  canonical anti-commutation relations :
\begin{align}
&\{ B (f) , B^{\ast}(g) \}  \, =  \,  (f, g)_{\ms{X}}  \, ,  \label{CAR1} \\
& \{ B (f) , B (g) \}  \, =   \{ B^{\ast} (f) , \,  B^{\ast }(g) \}  \,    = 0  , \label{CAR2}
\end{align}
where $ \{X, Y \} = XY + YX $.
Let $X$ be a self-adjoint operator on $\ms{X}$. Suppose that $X$ is bounded from below. Then the  second quantization  $\sqzf{X} $ is a self-adjoint  on $\Ff  (\ms{X}) $, which acts for the  vector 
 $\Psi = B^{\ast } (f_1) \cdots  B^{\ast } (f_n) \Omega_{\fer}$ as 
$  \sqzf{X}  \Psi $ $ = \sum\limits_{j=1}^n 
 B^{\ast } (f_1) \cdots  B^{\ast } (X f_j)  \cdots B^{\ast } (f_n) \Omega_{\fer}$. 
Let $f \in \ms{D}(X) $. Then it holds that  on  $\Ff^{\fin} (\ms{D} (X)) $,
\begin{align}
&[\sqzf{X} , B(f) ] \; = \; - \sqzf{Xf}  ,  \label{CRBf1} \\
&[\sqzf{X} , B^{\ast}(f) ] \; = \; \sqzf{Xf}  . \label{CRBf2}
\end{align}
 
$\;$ \\
Next we introduce  the boson Fock space.  The boson Fock space over a Hilbert space $\ms{Y}$ is defined by 
 $ \Fb (\ms{Y}) = \oplus_{n=0}^{\infty}
 (  \nstens \ms{Y} ) \, $, where $\nstens \ms{Y} $ denotes the $n$-fold  symmetric tensor product of $\ms{Y} $ with
$\; \tens_{\textrm{s}}^{0} \ms{Y} := \mbf{C}$.
The Fock vacuum  is defined by
$ \Omega_{\bos}  = \{   1, 0, 0, \cdots \} \, \in \, \Fb (\ms{Y}) $. 
The annihilation  operator is denoted by $A(f)$, $f \in \ms{Y}$, and  
the  creation operator by $ A^{\ast}(g)$, $ g \in \ms{Y}$.  
For a dense subspace $\ms{N} \subset \ms{Y} $, the finite particle subspace $\Fb^{\fin} (\ms{N} )$ is the linear hull of  vectors of the form $ \Psi =  A^{\ast } (f_1) \cdots  A^{\ast } (f_n) \Omega_{\bos}$, $ f_j \in \ms{N}$, $j=1, \cdots n$, $n \in \mbf{N} $.   
Creation operators and annihilation of bosonic field satisfy  canonical commutation relations on the finite particle subspace $\Fb^{\fin} (\ms{Y}) $: 
\begin{align}
&[ A (f) , A^{\ast}(g) ]  \, =  \,  (f, g)_{\ms{Y}}   ,    \label{CCR1} \\
& [ A (f) , A (g) ]  \, =   [ A^{\ast} (f) , A^{\ast }(g) ]  \,    = 0  , \label{CCR2}
\end{align}
 where $[ X, Y ] = XY -YX $.
Let $Y$ be a self-adjoint  operator on $\ms{Y} $. We assume that  $Y$ is non-negative. Then the  second quantization  $\sqzb{Y} $  is  self-adjoint on $\Fb (Y)$ which acts for the finite particle vector $ \Psi = A^{\ast } (f_1) \cdots  A^{\ast } (f_n) \Omega_{\bos}$
as 
$  \sqzb{Y}  \Psi $ $ = \sum\limits_{j=1}^n 
 A^{\ast } (f_1) \cdots  A^{\ast } (Y f_j)  \cdots A^{\ast } (f_n) \Omega_{\bos}$. 
     Let $f \in \ms{D}(Y) $. Then it follows that on  $\Fb^{\fin} (\ms{D} (Y)) $,
\begin{align}
&[\sqzb{Y} , A(f) ] \; = \; - \sqzb{Yf}   ,  \label{CRAf1}  \\
&[\sqzb{Y} , A^{\ast}(f) ] \; = \; \sqzb{Yf}    .  \label{CRAf2}
\end{align}
Let $f \; \in  \ms{D} (Y^{-1/2} )$. Then it follows that for all $\Psi \in \ms{D} (\sqzb{Y}^{1/2}) $,
\begin{align}
& \| A(f ) \Psi  \| \leq \| Y^{-1/2} f  \| \,
\|   \sqzb{X}^{1/2} \Psi \| , \label{bounda} \\
& \| A^{\ast}(f) \Psi  \| \leq \| Y^{-1/2} f \| \, \|   \sqzb{Y}^{1/2} \Psi \| + \| f  \| \| \Psi \| .
\label{boundad}
\end{align}

\subsection{Proof of Thorem\ref{main}}
\begin{lemma}   \label{strong-weak}
Let $\{ h_n \}_{n=1}^{\infty} $ be a sequence of $\ms{K} $ such that
w-$\lim\limits_{n \to \infty } h_n = 0$. 
Then for all $\Psi \in \Ff ( \ms{K} ) $,
\begin{align*}
& \textbf{(1)} \quad  \textrm{s-} \lim_{n \to \infty } B (h_n)  \Psi  = 0 ,   \qquad \\
& \textbf{(2)} \quad \textrm{w-} \lim_{n \to \infty } B^{\ast} (h_n)  \Psi  = 0 . 
\end{align*}
\end{lemma}
\textbf{(Proof)} \\
\textbf{(1)} Let $  \Psi =  B^{\ast } (f_1) \cdots  B^{\ast } (f_l) \Omega_{\fer}  \in \Ff^{\fin} (\ms{K})$ be a finite particle vector.
From canonical anti-commutation relations (\ref{CAR1}) and (\ref{CAR2}), it is seen  that 
\[
B(h_n) \Psi = \sum _{j=1}^l  (-1)^{j} (h_n , f_j   ) B^{\ast } (f_1) \cdots  \widehat{B^{\ast } (f_j )} \cdots B^{\ast } (f_l) \Omega_{\fer} , 
\]
where $\widehat{X} $ stands for omitting  the operator $X$. 
Since w-$\lim\limits_{n \to \infty } h_n = 0$, it follows that $\lim\limits_{n \to \infty} \| B(h_n) \Psi \| = 0 $.   Then we   see that  for 
all finite vector $\Psi \in \ms{F}_{\fer}^{\fin}( \ms{K} ) $, s-$\lim\limits_{n \to \infty } B(h_n) \Psi = 0 $. Note that $ \ms{F}_{\fer}^{\fin}( \ms{K} ) $ is dense in 
$ \Ff ( \ms{K} ) $ and  $\| B(h_n ) \| $ is uniformly bounded  with
 $ \| B(h_{n}) \| = \| h_n \| = 1 $, for all $n \in \mbf{N}$. Then we see that   s-$\lim\limits_{n \to \infty } B(h_n) \Psi = 0 $ for 
all  $\Psi \in \Ff ( \ms{K} )$. \\
\textbf{(2)} 
Let  $\Psi \in \Ff (\ms{K}) $.  From canonical anti-commutation relations (\ref{CAR1}) and (\ref{CAR2}), we see that for 
$  \Phi =  B^{\ast } (g_1) \cdots  B^{\ast } (g_l) \Omega_{\fer} \in \Ff^{\fin} (\ms{K})$, 
\[
(\Phi ,  B^{\ast}(h_n) \Psi )  = (B (h_n ) \Phi ,   \Psi ) 
= \sum _{j=1}^l  (-1)^{j} ( g_j , h_n  )  ( 
B^{\ast } (g_1) \cdots  \widehat{B^{\ast } (g_j)} \cdots B^{\ast } (g_l) \Omega_{\fer} , \Psi ).
\]
From this equality and    w-$\lim\limits_{n \to \infty } h_n = 0$,  we have 
$\lim\limits_{n \to \infty } (\Phi ,  B^{\ast}(h_n) \Psi ) = 0 $.
Note that $ \ms{F}_{\fer}^{\fin}( \ms{K} ) $ is dense in 
$ \Ff ( \ms{K} ) $ and $ \| B^{\ast}(h_{n}) \| = \| h_n \| = 1 $ for all $n \in \mbf{N}$. Hence we see that  for all  $\Phi \in \Ff ( \ms{K} )$, $\lim\limits_{n \to \infty } (\Phi ,  B^{\ast}(h_n) \Psi ) = 0 $.   $\blacksquare$

\begin{lemma}  \label{11/6.1}
It follows that for all $f \in \ms{D} (K)$, 
\begin{align*}
\textbf{(i)} \quad  [  H_0 , B (f) \tens I ]^0_{\ms{D} (H_0) } = - B (K f ) \tens I_{ \restr \ms{D}(H_0)} , \\
\textbf{(ii)}  \quad [  H_0 , B^{\ast} (f) \tens I ]^0_{\ms{D} (H_0) } =  B^{\ast}(K f ) \tens I_{ 
\restr \ms{D}(H_0)}  .
\end{align*}
\end{lemma}
\textbf{(Proof)} 
 Let $\Psi \in \Ff^{\fin} (\ms{D} (K)) \tens_{\textrm{alg}} \ms{D} (T) $, where  $\tens_{\textrm{alg}}$ denotes the algebraic tensor product, and  $\Phi \in \ms{D} (H_0 )$. From the commutation relation (\ref{CRBf1}),
\[
(H_0 \Phi , (B (f) \tens I  )  \Psi ) - ( (B^{\ast}(f) \tens I  )  \Phi , H_0 \Psi ) 
= -( \Phi , ( B (K f) \tens I ) \Psi ) .
\]
Note that $ \Ff^{\fin} (\ms{D} (K)) \tens_{\textrm{alg}} \ms{D} (T) $ is a core of $H_0$. 
Since $B(f)$ and $B^{\ast} (f) $ are bouded, we see that 
for all $ \Psi \in \ms{D} (H_0 )$,  $(H_0 \Phi , (B (f) \tens I  )  \Psi ) - ( (B^{\ast}(f) \tens I  )  \Phi , H_0 \Psi )  =-( \Phi , ( B (K f) \tens I ) \Psi ) $. Thus \textbf{(i)} holds. Similarly, 
we can prove  \textbf{(ii)}. $\blacksquare $

$\; $ \\
{\large \textbf{(Proof of Theorem \ref{main})} } \\
Let $ \lambda \in  \sigma_{\ess} (K) \backslash \{ 0  \}  $.
Then there exists a Weyl sequence $\{h_n \}_{n=1}^{\infty}$ for $K$ and $\lambda $, i.e., 
$h_n \in \ms{D} (K)$ and $\| h_n \| =1 $ for all $n \in \mbf{N}$, w-$\lim\limits_{n \to \infty} h_n = 0 $, s-$\lim\limits_{n \to \infty} (K -\lambda ) h_n= 0 $.  
By  this sequence, we construct  a  Weyl sequence
for $H$ and $E_0 (H) + \lambda $ as follows.
 Let us set   
\begin{equation}
\Psi_{n , \epsilon } \; =  \; \left( \left(  B(h_n ) + B^{\ast} (h_n) \right) \tens I \right) \Xi_\epsilon   ,
\end{equation}
where $ \Xi_\epsilon  \in \text{Ran} (E_{H} ( [ 0, \epsilon )) )$, $\| \Xi_\epsilon \| =1$,
$ 0 < \epsilon \leq 1$. Here $E_{H} $ denotes the spectral projection of $H$.
From canonical anti-commutation relations (\ref{CAR1}) and (\ref{CAR2}),
\begin{align}
\| \Psi_{n , \epsilon }  \|^2 &   = (  \Xi_\epsilon  , \left( (B(h_n )^2  +  \left\{  B(h_n ) , B^{\ast} (h_n) \right\}  + B^{\ast}(h_n )^2) \tens I \right)
\Xi_{\epsilon} )  \notag  \\  
& = \| h_n \|^2   \, \| \Xi_{\epsilon}  \|^2 . 
\end{align}
Since $\| h_n \| =1 $ and $\| \Xi_{\epsilon}\|=1$,   we see that  $\| \Psi_{n , \epsilon }  \| = 1$ for all $n \geq 1$ and $ 0 < \epsilon  \leq 1 $.  
From Lemma \ref{11/6.1} and the assumption \textbf{(A.2)}, it holds that for all
$\Phi  , \Theta \in \ms{D} (H ) $ and for all $h \in \ms{D} (K)$, 
\begin{align*}
& ( H \Phi,     ((B (h)  + B^\ast (h) ) \tens I ) \Theta )  
- ( ((B (h)  + B^\ast (h) ) \tens I ) \Phi , H  \Theta  ) 
\\ 
&\quad = \left( \Phi, \left(  ( B^{\ast} (Kh) - B (Kh)  )\tens I   + [\HI , B^{\ast} (h) \tens I ]^{0}_{\ms{D}(H)}  +  [\HI , B(h) \tens I  ]^{0}_{\ms{D}(H)} \right) \Theta \right)  .
\end{align*}
Hence we see that  for all $\Phi \in \ms{D} (H)$, 
\begin{align}
(H \Phi , \Psi_{n , \epsilon } )  &=  (\Phi, ( (B (h_n) + B^{\ast} (h_n ) ) \tens I )  H \Xi_\epsilon )
+  (\Phi, ( (B^{\ast} (Kh_n )   - B (K h_n ) ) \tens I ) \Xi_\epsilon )  \notag  \\
& \quad  + (\Phi, [ \HI   , (B^{\ast} (h_n)  \tens I ) ]^0_{\ms{D}(H)}  \Xi_\epsilon ) + (\Phi, [ \HI   , (B (h_n)  \tens I ) ]^0_{\ms{D}(H)} \Xi_\epsilon  ) .
\end{align}
Then  $ \Psi_{n , \epsilon } \in \ms{D} (H^\ast )$,  and hence, 
$ \Psi \in \ms{D} (H) $, since $H$ is self-adjoint. Then we have
\begin{align}
 H \Psi_{n , \epsilon }  =& 
 ( ( B (h_n) + B^{\ast} (h_n ) ) \tens I  )H \Xi_\epsilon 
+ ( (B^{\ast} (Kh_n )   - B (K h_n ) ) \tens I )  
 \Xi_\epsilon  \notag \\
 & +  [\HI   ,  B^{\ast} (h_n)  \tens I ]^0_{\ms{D}(H)} \Xi_\epsilon  
+  [\HI   , B(h_n)  \tens I ]^0_{\ms{D}(H)} \Xi_\epsilon  . \label{12/27-1}
\end{align}
From (\ref{12/27-1}),   
\begin{align}
&  \|  \left(  H -   (\lambda + E_0 (H)) \right)  \Psi_{n , \epsilon} \|   \notag \\
 &\quad    \leq \|  ( (B (h_n) + B^{\ast} (h_n ) ) \tens I ) (H - E_0 (H)) \Xi_\epsilon \|
+ \| ( B^{\ast} (K h_n  - \lambda h_n ) \tens I ) \Xi_{\epsilon}  \|   \notag  \\
&\quad \quad + \| B (K h_n  + \lambda h_n ) \Xi_{\epsilon}  \| 
+  \|  [H   , B (h_n)  \tens I ]^0_{\ms{D}(H)}  \Xi_\epsilon  \| + 
 \|  [H   , B^{\ast}(h_n)  \tens I ]^0_{\ms{D}(H)}  \Xi_\epsilon    \|  . \label{11/6.abc}
\end{align}
We see that 
\begin{equation}
\| ( (B (h_n ) +  B^{\ast } (h_n)  ) \tens I )( H-E_0 (H) )  \Xi_\epsilon \|  
\leq 2 \| h_n \| \|  ( H-E_0 (H) )  \Xi_\epsilon \|   \leq 2 \epsilon  , \label{11/6.a}
\end{equation}
and 
\begin{equation}
\| ( B^{\ast} (K h_n  - \lambda h_n ) \tens I ) \Xi_{\epsilon}  \|  \leq  \| (K   - \lambda )h_n\| .
\label{11/6.b}
\end{equation}
Since $  B (K + \lambda ) h_n =  B(K-  \lambda ) h_n    + 2 \lambda   B (h_n)   $, we also see that
\begin{equation}
\| (B (  (K+ \lambda ) h_n ) \tens I )  \Xi_\epsilon\| \leq 
\|  (K- \lambda ) h_n \|    + 2 | \lambda | \|   (B (h_n) \tens I) \Xi_\epsilon   \| . \label{11/6.c}
 \end{equation}
By applying  (\ref{11/6.a}), (\ref{11/6.b}) and (\ref{11/6.c}) to (\ref{11/6.abc}), 
\begin{equation}
 \|  \left(  H -   (\lambda + E_0 (H)) \right)  \Psi_{n , \epsilon} \|  
\leq 2 \epsilon + 2 \|   (K- \lambda ) h_n  \|  + 2 | \lambda | 
 \|   (B (h_n) \tens I )\Xi_\epsilon \| .  \label{11/6.2}
\end{equation}
Since s-$ \lim\limits_{n \to \infty} (K- \lambda ) h_n = 0 $ and  s-$ \lim\limits_{n \to \infty}  (B (h_n) \tens I ) \Xi_\epsilon = 0$ by Lemma \ref{strong-weak}, 
(\ref{11/6.2}) yields that  
\begin{equation}
\lim_{\epsilon \to 0} \limsup_{n \to \infty } \|  \left(  H -   (\lambda + E_0 (H)) \right)  \Psi_{n , \epsilon} \|  = 0 . \notag 
\end{equation}
Then, we can take a subsequence $\{\Psi_{n_j  , \epsilon_j } \}_{j=1}^{\infty}$ satisfying 
$\lim\limits_{j \to \infty } \|  \left(  H -   (\lambda + E_0 (H)) \right)  \Psi_{n_j , {\epsilon}_{j}} \|  $ $ = 0$. In addition, from the definition of $\Psi_{n, j}$ and Lemma \ref{strong-weak}, we see that  w-$\lim\limits_{j \to \infty}\Psi_{n_j  , \epsilon_j } =0$. Thus $\{\Psi_{n_j  , \epsilon_j } \}_{j=1}^{\infty}$ is a Weyl sequence for $E_0 ( H )  + \lambda  $.
Then  Weyl's criterion  (\cite{HiSi96}) says that  $  E_0 ( H )  + \lambda      \in  \sigma_{\ess} (H) $. Since $\sigma_{\ess} (H)$ is closed in $\mbf{R}$,  the  proof is obtained.  $\blacksquare$

\section{Application}
We consider an application of Theorem \ref{main} to  the system of a Dirac field interacting with a Klein-Gordon  field. The state space of Dirac field and Klein-Gordon field are given by $ \ms{H}_{\textrm{D}} = \Ff (\stateDirac )$ and $\ms{H}_{\textrm{KG}} = \Fb (\stateKG ) $, respectively. The state space is 
defined by 
 \[
 \ms{H}_{\textrm{Y}}  \; = \ms{H}_{\textrm{D}}  \tens \ms{H}_{\textrm{KG}} .
 \]
The free Hamiltonians  of the Dirac field and the Klein-Gordon field are given by  $\HDirac =  \sqzf{\omega_{M}} $ with  $\omega_{M} (\mbf{p}) = \sqrt{ \mbf{p}^2 + M^2} $,  $M>0 $, and $ \HKG \; = \; \sqzb{\omega_{m}} $ with $\omega_{m} (\mbf{k}) = \sqrt{ \mbf{k}^2 + m^2} $, $m>0 $, respectively. The total Hamiltonian is defined by 
\[
H_{\kappa} \; = \;   \HDirac \tens I + I  \tens  \HKG   + \kappa \HI  , \qquad  \kappa \in \mbf{R} , 
\]
where $\HI$ is  the symmetric operator on $\ms{H}_{\Yukawa}$ such that
 for all $\Phi \in \ms{H}_{\Yukawa} $, and for all
$ \Psi \in \ms{D} (H_0)$, 
\begin{equation}
 (\Phi ,  \HI \Psi  ) \; = \;
\int_{\Rthree} \cutoff{I}{x} \left( \Phi , \,
 \overline{\psi (\mbf{x}) } \psi (\mbf{x} ) \tens \phi (\mbf{x})
 \Psi  \right) \dx . 
 \end{equation}
Here $\psi (\mbf{x} )   $ and $\phi (\mbf{x})$ are field operators of the Dirac field and 
 Klein-Gordon field, respectively, and
$  \overline{\psi (\mbf{x}) } \, = \, \psi^{\ast}(\mbf{x} ) \beta$. Dirac matrices $\alpha^{j} $, $j=1,2,3$, and  $\beta$   are the $4\times 4$ hermitian matrices satisfying  anti-commutation relations $
 \{ \alpha^{j} , \alpha^{l} \} = 2 \delta_{j,l} , \;
\{ \alpha_{j} , \beta \} = 0 $ and $   \; \beta^{2}=I  $. 
The definitions of  $\psi (\mbf{x} )   $ and $\phi (\mbf{x})$ are as follows.  First we consider the Dirac field's operator.
 Let $B(\xi )$, $ \; \xi = (\xi_{1}, \cdots , \xi_{4} )  \in \ms{H}_{\Dirac}$, and
$\; B^{\ast}(\eta)$, $ \; \eta = (\eta_{1}, \cdots , \eta_{4} )  \in \ms{H}_{\Dirac}$, be the
 annihilation operator and the creation operator on  $\ms{H}_{\Dirac}$, respectively.
  Let us set for $f \in \eltwo$,
\begin{align*}
&b^{\ast}_{1/2}(f) = B^{\ast}( (f,0,0,0) ), \quad \quad
b^{\ast}_{-1/2}(f) = B^{\ast}((0,f,0,0) ), \\
&d^{\ast}_{1/2}(f) = B^{\ast}((0,0,f,0) ), \quad \quad
d^{\ast}_{-1/2}(f) = B^{\ast}((0,0,0,f) ) .
\end{align*}
Then they satisfy canonical anti-commutation relations below: 
 \begin{align}
&\{ b_{s} (f) , b_{\tau}^{\ast}(g) \} =
\{ d_{s} (f) , d_{\tau}^{\ast}(g) \} =
\delta_{s, \tau} (f, g)_{\eltwo}   ,    \label{CARDirac1} \\
&\{ b_{s} (f) , b_{\tau} (g) \} =
\{ d_{s} (f) , d_{\tau} (g) \} = \{ b_{s} (f) , d_{\tau} (g) \} =
\{ b_{s} (f) , d_{\tau}^{\ast}(g) \} = 0  . \label{CARDirac2}
\end{align}
Let  $ u_{s}(\mbf{p} ) =
(u_{s}^{l} (\mbf{p}) )_{l=1}^{4} \; $ and  $v_{s}(\mbf{p} ) = (v_{s}^{l} (\mbf{p}) )_{l=1}^{4}  \; $ be  spinors   with spin $s =  \pm 1/2 $, which are the positive and negative energy part of the Fourier transformed Dirac operator $ h_{\Dirac} (\mbf{p}) = \mbf{\alpha} \mbf{\cdot} 
  \mbf{p}   + \beta M $, respectively. 
The filed operator  $\psi(\mbf{x}) = ( \psi_{l}(\mbf{x}) )_{l=1}^4$ of the Dirac field is defined by 
\[
\psi_{l}(\mbf{x}) = \sum_{s=\pm 1/2}( b_{s} (f_{s,\mbf{x}}^{l} ) +   d^{\ast}_{s} (g_{s, \mbf{x}}^{l} )) ,
 \quad   l= 1, \cdots  ,4 ,
\]
where $\;   f_{s,\mbf{x}}^{l} (\mbf{p}) = f_{s}^{l} (\mbf{p})
e^{-i \mbf{p} \cdot \mbf{x}} \; $ with $f_{s}^{l} (\mbf{p}) =
\frac{\chi_{\Dir}(\mbf{p} ) u_{s}^{l} (\mbf{p}) }{ \sqrt{ (2 \pi ) ^3  \omega_M (\mbf{p})}}$ 
 and $ \;  g_{s,\mbf{x}}^{l} (\mbf{p}) = g_{s}^{l} (\mbf{p})
e^{-i \mbf{p} \cdot \mbf{x}} $  with $g_{s}^{l} (\mbf{p})=
\frac{\chi_{\Dir} (\mbf{p})    v_{s}^{l} (-\mbf{p})}
{\sqrt{ (2 \pi ) ^3  \omega_M (\mbf{p} )}}  $. Here  $\chi_{\Dir}$ denotes  an  ultraviolet cutoff.
Next we define the Klein-Gordon filed's operators. Let $a(f)$, $f \in \ms{H}_{\KG} $, and $a^{\ast}(g), g \in \ms{H}_{\KG}  $, be the annihilation operator and the 
creation operator, respectively. Then it is seen that on $\Fb^{\fin} (\ms{H}_{\KG} )$, 
\begin{equation}
[ a (f) , a^{\ast} (g) ] = 0 \, \quad  \; [a(f) , a(g) ] = [a^{\ast} (f) , a^{\ast} (g) ] = 0. \label{CCRKG}
\end{equation}
The field operator  of the Klein-Gordon field is defined by
\[
 \phi (\mbf{x}  ) \; = \; \frac{1}{\sqrt{2}}
 \left( \frac{}{} a (h_{\mbf{x}}) \; + \;  a^{\ast} (h_{\mbf{x}}) \right),
\]
where $h_{\mbf{x}}  (\mbf{k})  = h (\mbf{k}) e^{i\mbf{k} \cdot \mbf{x}} \; $ with
$ \; h (\mbf{k}) \; = \; \frac{\cutoff{KG}{k} }{\sqrt{(2 \pi )^3 \omega_m  (\mbf{k}) }}$,
and $\chi_{\KG}$ is an  ultraviolet cutoff.

$\;$ \\
We suppose the following conditions.
\begin{quote}
\textbf{(Y.1) (Ultraviolet cutoff)}
\[
\int_{\Rthree}  \left| \chi_{\Dir} (\mbf{p})
  u_{s}^{l} (\mbf{p}) \right|^2  d \mbf{p}  \; <  \; \infty  , 
 \left| \chi_{\Dir} (\mbf{p}) v_{s}^{l} (- \mbf{p}) \right|^2  d \mbf{p}  \; <  \; \infty , 
\int_{\Rthree}
  \left| \cutoff{KG}{k} \right|^2   \dk \; < \; \infty  .
 \]
\textbf{(Y.2)}  \textbf{(Spatial cutoff)}
 $ \int_{\Rthree} | \chi_{\I} (\mbf{x}) | \dx < \infty$  .
\end{quote}

$\;$ \\
From the boundedness (\ref{boundBBast}), (\ref{bounda}), (\ref{boundad}) and $ \overline{\psi (\mbf{x}) } \psi (\mbf{x} ) = \sum_{l, l'} \beta_{l,l'} \psi_l^{\ast} (\mbf{x} )
 \psi_{l'} (\mbf{x} ) $, we see by
 \textbf{(Y.1)} that 
\begin{align}
& \sup_{\mbf{x} \in \Rthree} \| \overline{\psi (\mbf{x})}  \psi (\mbf{x}) \|  \leq 
 \sum_{l,l' =1}^4  \, \sum_{s, s'=\pm 1/2} |\beta_{l , l'} |  (  \|     f_{s}^{l}  \| + \| g_{s}^{l}  \| )   (  \|     f_{s'}^{l}  \| + \| g_{s'}^{l}  \| )  , \label{boundpsi}   \\
& \sup_{\mbf{x} \in \Rthree} \| \phi (\mbf{x}) \Psi  \|
\leq \sqrt{2}  \|  \frac{h}{\sqrt{\omega_{m}}} \|  \, \|  \HKG^{1/2} \Psi  \|  + 
 \frac{1}{\sqrt{2}}\| h \| \|  \Psi\|  . 
\label{boundphi}
 \end{align}
 Then from (\ref{boundpsi}), (\ref{boundphi}), \textbf{(Y.2)} and $ \| \HKG^{1/2} \Phi  \| \leq  \epsilon \| \HKG \Phi \| + \frac{1}{2 \epsilon} \| \Phi \|$, it holds that for all $\Psi \in \ms{D}(\Hzero )$,
\begin{align*}
\| \HI \Psi \|   & \leq  
\epsilon \sqrt{2} \| \chi_{\I} \|_{L^1 }  \| \frac{h}{\omega_{\KG}} \| \| \Hzero  \Psi \|   
  + \| \chi_{\I} \|_{L^1 }  \left(  \frac{1}{ \sqrt{2} \epsilon } \| \frac{h}{\omega_{m}} \|  
+ \frac{1}{\sqrt{2}}  \| h \|  \right. \\
&  \qquad \quad   \; + \left. \sum\limits_{l,l' =1}^4  \, \sum\limits_{s, s'=\pm 1/2} |\beta_{l , l'} |  (  \|     f_{s}^{l}  \| 
 + \| g_{s}^{l}  \| )   (  \|     f_{s'}^{l}  \| + \| g_{s'}^{l}  \| ) \right) \| \Psi  \| .
\end{align*}
Thus $\HI $ is relatively bounded with respect to $H_0 $ with the infinitely small bound, and then the Kato-Rellich theorem shows that  $H_{\kappa} $ is self-adjoint and essentially self-adjoint on any core of $\Hzero $. In particular $H_{\kappa} $ is essentially self-adjoint on
$\Ff^{\fin} (\ms{D}(\omega_{\Dirac}) ) \tens_{\textrm{alg} } \Fb^{\fin} (\ms{D} (\omega_{\KG}))  $,
where $\tens_{\textrm{alg} }$ stands for the algebraic tensor product.

$\;$ \\
Let $ \nu  = \min \{ m, M  \}$. From  Theorem\ref{main}, the next assertion  follows.  
\begin{theorem} \label{Application} $\;$ \\
Assume \textbf{(Y.1)} and \textbf{(Y.2)}. Then  
$   [ E_0 (H_{\kappa } ) + \nu , \, \infty     ) \subset  \sigma_{\ess} (H_{\kappa})$
for all $ \kappa \in \mbf{R}$.  
\end{theorem}

$\;$ \\
Before proving   Theorem \ref{Application}, we explain a  result of the essential spectrum of the Yukawa model and state a corollary. In \cite{Ta11} the following theorem has been  proven: 
\begin{quote}
\textbf{Theorem A (\cite{Ta11} ; Theorem 2.1)}  \\
Assume \textbf{(Y.1)},  \textbf{(Y.2)} and $ \int_{\Rthree} |\mbf{x}| \, | \chi_{\I} (\mbf{x}) | \dx < \infty$. Then it follows that for all $ \kappa \in \mbf{R}$, 
$  \sigma_{\ess} (H_{\kappa}) \subset [ E_0 (H_{\kappa } ) + \nu  , \, \infty  )  $.
\end{quote}

$\;$ \\
From Theorem A and Theorem \ref{Application}, the next corollary  follows. 
\begin{corollary} (\textbf{HVZ theorem for the Yukawa model}) \\
Assume \textbf{(Y.1)},  \textbf{(Y.2)} and $ \int_{\Rthree} |\mbf{x}| \, | \chi_{\I} (\mbf{x}) | \dx < \infty$. Then 
$   [ E_0 (H_{\kappa } ) + \nu , \, \infty     ) =  \sigma_{\ess} (H_{\kappa})$
 for all $ \kappa \in \mbf{R}$.
\end{corollary}

$\, $ \\
To prove  Theorem \ref{Application}, we prepare for some lemmas. 
\begin{lemma}  \label{Generalized}
Let  $A$ and $B$ be    self-adjoint  operators on Hilbert spaces $\ms{X}$ and 
$\ms{Y}$, respectively. Assume that $A$ and $B$ are non-negative. 
Let 
$X ( \mbf{x} )  $, $ \mbf{x} \in \Rd$, be a linear  operator on $\Ff (\ms{X}) $, and $Y(\mbf{x})$, $ \mbf{y} \in \Rd$, a linear  operator on $\Fb (\ms{Y})$, which satisfy 
\begin{align}
&\quad \sup_{\mbf{x} \in \Rd }\|    X ( \mbf{x} )   \Psi    \|   \leq  \textrm{const.}
(  \|  \sqzf{A}^{1/2}  \Psi \|  +    \|  \Psi \| )  ,      \quad \Psi \in \ms{D} ( \sqzf{A}^{1/2} ) ,  \\
&\quad  \sup_{\mbf{x} \in \Rd}  \|    Y ( \mbf{x} )   \Psi    \|   \leq
\textrm{const.} (  \|  \sqzb{B}^{1/2} \Psi  \| +   \|  \Psi \| ), \quad  \Psi \in \ms{D} ( \sqzb{B}^{1/2} ) ,
\end{align}
respectively.  
Then  there exists a linear operator $Z $ on $\Ff (\ms{X})  \tens \Fb (\ms{Y})$  such that
 $\ms{D} ( \sqzf{A}   \tens I ) \cap \ms{D} (I \tens \sqzb{B} )  \subset \ms{D} (Z)$  and 
 for all $\Phi \in \Ff (\ms{X})  \tens \Fb (\ms{Y})$ and for all $\Psi \in \ms{D} ( \sqzf{A}   \tens I ) \cap \ms{D} (I \tens \sqzb{B} )$,
\[
 (\Phi , Z \Psi) = 
 \; \int_{\Rd}  g (\mbf{x})  ( \Phi ,   (  X ( \mbf{x} )  \tens  Y ( \mbf{x} )   \Psi ) \dx ,
\]
where  $g $ is the Borel function satisfying  $\| g\|_{L^1} < \infty$. 
\end{lemma}
\textbf{(Proof)}  Let  $\ell $ be a linear  functional on $ (\Ff (\ms{X} ) \tens \Fb (\ms{Y}) ) \times  (\Ff ( \ms{D} (A) ) \tens \Fb ( \ms{D} (B) ) )$ defined by
\begin{equation}
\ell (\Phi , \Psi ) \; = \; \int_{\Rd}  g (\mbf{x})  ( \Phi ,   (  X ( \mbf{x} )  \tens  Y ( \mbf{x} )   \Psi ) \dx . 
\end{equation} 
It is seen that $ | \ell (\Phi , \Psi )  |  \leq $ $c \|g \|_{L^1} $ $\| \Phi \| 
( \| ( \sqzf{A}  \tens I    + I \tens \sqzb{B}   \Psi \| + \| \Psi\| )$ with some constant $c>0$. From the Riez representation theorem, the assertion holds.   $\blacksquare$

$\;$ \\
From  canonical anti-commutation relations (\ref{CARDirac1}) and (\ref{CARDirac2}),  it is proven that commutation relations below
follow  in (\cite{Ta09} ; Lemma 3.1)  : 
\begin{align}
&[ \psi_{l}^{\ast}(\mbf{x}) \, \psi_{l'}(\mbf{x})  , \; b_{s }(\xi  )  ] = -( \xi, f_{s, \mbf{x}}^{l}  )
 \psi_{l'}(\mbf{x})   ,   \label{ppbs} \\
&[ \psi_{l}^{\ast}(\mbf{x}) \, \psi_{l'}(\mbf{x})  , \; d_{s}(\xi )  ] = ( \xi , g_{s , \mbf{x}}^{l'}  )  
\psi_{l}^{\ast}(\mbf{x}) .\label{ppds}
\end{align}

$\;$ \\
From (\ref{ppbs}), (\ref{ppds}),   and  $[X, Y]^{\ast} = - [X^\ast , Y^\ast  ]$, it is seen that 
\begin{align}
&[ \psi_{r}(\mbf{x}) \, \psi_{r'}^{\ast}(\mbf{x})  , \; b^{\ast}_{s }(  \eta)  ] = 
 (  f_{s, \mbf{x}}^{r'}  , \eta  )   \psi_{r}^{\ast}(\mbf{x})   , \label{ppbs*}  \\
&[ \psi_{r} (\mbf{x}) \, \psi_{r'}^{\ast}(\mbf{x})  , \; d_{s}^{\ast} (\eta)  ] =
 - (  g_{s , \mbf{x}}^{r}  , \eta )  \psi_{r'} (\mbf{x}) . \label{ppds*}
\end{align}

$\;$ \\
Here note that
\begin{equation}
\sup\limits_{\mbf{x} \in \Rthree} \|   \psi_l (\mbf{x}) \|  \leq \sum\limits_{s =\pm 1/2}  (  \|     f_{s}^{l}  \| + \| g_{s}^{l}  \| ) . \label{11/19.1}
\end{equation}
From  commutation relations (\ref{ppbs}) - (\ref{ppds*}) and bounds   (\ref{boundphi}) and  (\ref{11/19.1}), Lemma \ref{Generalized} yields that  the next lemma follows.
\begin{lemma} \label{CommutationCAR}
$\;$ \\
\textbf{(1)}There exist $ [  \HI ,  b_{s} (\xi )  \tens I ]^0_{\ms{D}(H)}  $ and 
$ [  \HI ,  d_{s} (\xi )  \tens I ]^0_{\ms{D}(H)}  $
such that  for all $\Phi \in \HYukawa$ and for all $ \Psi \in \ms{D} (\Hzero )$,
\begin{align}
&  ( \Phi , [  \HI ,  b_{s} (\xi )  \tens I ]^0_{\ms{D}(H)} \Psi ) \; = \; 
- \sum_{l, l'} \beta_{l, l'}\int_{\Rthree}  \chi_{\I} (\mbf{x})  ( \xi, f_{s, \mbf{x}}^{l}  )
  \left( \Phi , (  \psi_{l'} (\mbf{x} )  \tens  \phi (\mbf{x})  ) \Psi  \right) \dx   ,  \label{HIbs} \\
& ( \Phi , [  \HI ,  d_{s} (\xi )  \tens I ]^0_{\ms{D}(H)} \Psi ) \; = \; 
\sum_{l, l'} \beta_{l, l'}\int_{\Rthree}  \chi_{\I} (\mbf{x})   ( \xi , g_{s , \mbf{x}}^{l'}  )
 \left( \Phi ,  \psi_{l}^{\ast} (\mbf{x}) \tens   \phi (\mbf{x}) ) \Psi \right)  \dx  . \label{HIds}
\end{align}
\textbf{(2)} There exist $ [  \HI ,  b_{s}^{\ast} (\eta )  \tens I  ]^0_{\ms{D}(H)}  $ and 
$ [  \HI ,  d_{s}^{\ast} (\eta )  \tens I ]^0_{\ms{D}(H)}  $
such that for all $\Phi \in \HYukawa$ and for all $ \Psi \in \ms{D} (\Hzero )$
\begin{align}
& ( \Phi , [  \HI ,  b^{\ast}_{s} (\eta )  \tens I ]^0_{\ms{D}(H)} \Psi ) \; = \; 
\sum_{r, r'} \beta_{r, r'} \int_{\Rthree}  \chi_{\I} (\mbf{x})  (  f_{s, \mbf{x}}^{r'}  , \eta  )
 \left( \Phi , (   \psi_{r}^{\ast} (\mbf{x} )  \tens  \phi (\mbf{x})  ) \Psi \right) \dx  , 
\label{HIbs*}  \\
& ( \Phi , [  \HI ,  d^{\ast}_{s} (\eta )  \tens I ]^0_{\ms{D}(H)} \Psi ) \; = \; 
-  \sum_{r, r'} \beta_{r, r'} \int_{\Rthree}  \chi_{\I} (\mbf{x})  ( g_{s , \mbf{x}}^{r}  , \eta )
  \left( \Phi ( \psi_{r'} (\mbf{x})  \tens \phi (\mbf{x})   ) \Psi  \right)    \dx  .   \label{HIds*}
\end{align}
\end{lemma}

$\; $ \\
Similarly,  the next  proposition follows from  (\ref{CCRKG}) and  (\ref{boundpsi}) to Lemma \ref{Generalized}.
\begin{lemma} \label{CommutationCCR}
There exits $ [  \HI , I \tens a^{\ast} (\zeta) ]^0_{\ms{D}(H)}  $ satisfying for all $\Phi \in H_{\Yukawa}$ and for all $ \Psi \in \ms{D} (\Hzero )$,
\[ 
(\Phi, [  \HI , I \tens a^{\ast} (\zeta) ]^0_{\ms{D}(H)}   \Psi ) \; = \; 
\int_{\Rthree}  \chi_{\I} (\mbf{x})  (f_{\mbf{x} }  , \zeta ) 
(\Phi ,  (  \overline{\psi (\mbf{x}) } \psi (\mbf{x} )  \tens  I  ) \Psi ) \dx .
\]
\end{lemma}

$\; $ \\
{\large \textbf{(Proof of Theorem \ref{Application})}}  \\
Let us apply to  Theorem \ref{main} to the Yukawa model. Since
$H_{\kappa}$ is self-adjoint on $\ms{D}(H_{0})$ and  bounded from below, 
\textbf{(A.1)} is satisfied. Let us check  \textbf{(A.2)}. 
Let  $ \{ h_n \}_{n=1}^\infty $ be a sequence such that
 $h_n \in \ms{D}(\omega_{M}) , \, \|h_n \| = 1 $ for all $n \in \mbf{N}$, and w-$\lim\limits_{n \to \infty} h_n = 0$. 
From (\ref{HIbs}), it is seen that  for all $\Psi \in \ms{D} (H)$, 
\begin{equation}
\|  [  \HI ,  b_{s} (h_n )  \tens I ]^0_{\ms{D}(H)} \Psi \|  \; \leq    \sum_{l, l'} | \beta_{l, l'} | \int_{\Rthree}  \chi_{\I} (\mbf{x})  | ( h_n , f_{s, \mbf{x}}^{l}  )| \, 
\|   (  \psi_{l'} (\mbf{x} )  \tens  \phi (\mbf{x})  ) \Psi \|  \dx  \\  . \label{11/9.1}
\end{equation}
By  (\ref{11/19.1}) and (\ref{boundphi}), we have 
$  \sup\limits_{\mbf{x}\in \Rthree}  \|  (  \psi_{l'} (\mbf{x} )  \tens  \phi (\mbf{x})  ) \Psi \| < \infty  $. We also see that $\int_{\Rthree} |\chi_{\I} (\mbf{x}) |d \mbf{x} < \infty$ by \textbf{(A.1)} and $ | ( h_n , f_{s, \mbf{x}}^{l}  )|  \leq  \|f_{s} \|  $.  
 Then from w-$\lim\limits_{n \to \infty} h_n =0 $ and  the Lebesgue dominated convergence theorem,  (\ref{11/9.1}) yields that  $\lim\limits_{n \to \infty} \|  [  \HI ,  b_{s} (h_n )  \tens I ]^0_{\ms{D}(H)} \Psi \| =0$, $s= \pm \frac{1}{2}$. Similarly we can prove that  $\lim\limits_{n \to \infty} \|  [  \HI ,  d_s (h_n )  \tens I ]^0_{\ms{D}(H)} \Psi \| =0$, $s= \pm \frac{1}{2}$ . Then the condition (1)
in \textbf{(A.2)} is satisfied. In addition, we can also prove   that for  $s= \pm \frac{1}{2}$,
$\lim\limits_{n \to \infty} \|  [  \HI ,  b_s^{\ast} (h_n )  \tens I ]^0_{\ms{D}(H)} \Psi \| =0$
 and $\lim\limits_{n \to \infty} \|  [  \HI ,  d_s^{\ast} (h_n )  \tens I ]^0_{\ms{D}(H)} \Psi \| =0$, and then the condition (2)
in \textbf{(A.2)} is satisfied.  
 Hence from Theorem \ref{main}, it follows that  
$   [ E_0 (H_{\kappa } ) + M  , \, \infty   ) \subset  \sigma_{\ess} (H_{\kappa})$. 
Next  we show $   [ E_0 (H_{\kappa } ) + m   , \,  \infty ) \subset  \sigma_{\ess} (H_{\kappa})$, and then  Theorem \ref{Application} is proven. Here note that $\omega_{m}^{-r} $ is bounded for all $r>0$,
 since $\omega_{m} > 0 $.
Let $ \{f_n \}_{n=1}^\infty $  be a sequence of $\ms{H}_{m}  $ such that
 $f_n \in \ms{D}(\omega_{m}) , \, \|f_n \| = 1 $ for all $n \in \mbf{N}$, and w-$\lim\limits_{n \to \infty} h_n = 0$.  From  Lemma \ref{CommutationCCR}, it is seen that for all $\Psi \in \ms{D} (H)$,
\begin{equation}
  \| [  \HI , I \tens a^{\ast} (f_n ) ]^0_{\ms{D}(H)}   \Psi \| \; 
  \leq \;  \int_{\Rthree}  | \chi_{\I} (\mbf{x}) |  \, | (f_{\mbf{x} }  , f_n ) |  \| (  \overline{\psi (\mbf{x}) } \psi (\mbf{x} )  \tens  I  ) \Psi \|  d \mbf{x}\; .  \label{11/15.a}
\end{equation}
We see that  $\chi_{I} \in L^1$, $|  (h_{\mbf{x} }  , f_n ) | \leq \| f\| $,  w-$\lim\limits_{n \to \infty} f_n = 0$ and  $\sup\limits_{\mbf{x} \in \Rthree} \| \overline{\psi (\mbf{x}) } \psi (\mbf{x} )  \|  < \infty $. Then 
we have $\lim\limits_{n \to \infty} \| [  \HI , I \tens a^{\ast} (f_n ) ]^0_{\ms{D}(H)}   \Psi \| \;  = 0 $ 
by (\ref{11/15.a}) and the Lebesgue dominated convergence theorem. 
 Then  \textbf{(S.2)} of Theorem I  in Appendix is satisfied, and hence $   [ E_0 (H_{\kappa } ) + m   , \,  \infty ) \subset  \sigma_{\ess} (H_{\kappa})$. 
 $\blacksquare $

$\;$ \\
$\;$ \\
{\large \textbf{Appendix}  (\cite{Arai00};Theorem 1.2) }  \\
$\; $ \\
Let
\[
\ms{H} \; =  \ms{R}  \tens    \Fb (\ms{S} ) 
\]
where $\ms{R}$ is a 
Hilbert pace and $\Fb ( \ms{S} )$  the boson Fock space over a Hilbert space $\ms{S}$.   Let  $R$ be a self-adjoint operator on $\ms{R}$ and $S$  a self-adjoint operator on $\ms{S}$  with ker $S = \{ 0 \} $. Additionally, we assume that  $R$ is bounded from below and $S$ is non-negative. Let
 \[
H \; = \; R  \tens I \; + \; I \tens  \sqzb{S}  \; +  \; \HI   ,
\]
where $\HI$ is a symmetric operator on $\ms{H}$. Let $\Hzero = R  \tens I + I \tens  \sqzb{S}$. We suppose  following conditions : 
\begin{quote}
\textbf{(S.1)} $H$ is self-adjoint on $\ms{D} (H) = \ms{D} (H_0) \cap \ms{D} (\HI )$ and bounded from below. \\
\textbf{(S.2)} For all $h \in \ms{D} (S) \cap \ms{D} (S^{-1/2})$, the weak commutator $[\HI  , I \tens A^{\ast} (h) ]^0_{\ms{D}(H)} $ exists. Moreover, for all sequences $ \{ h_n \}_{n=1}^\infty $ of $\ms{D} (S)  \cap \ms{D} (S^{-1/2}) $ such that $\| f_n \| = 1 $, $n \geq 1$, and w-$\lim\limits_{n \to \infty} h_n =0$,    it follows that for all $\Psi \in \ms{D} (H)$,
\[
\textrm{s-}\lim_{n \to \infty}[\HI  , I \tens A^\ast (h_n)   ]^0_{\ms{D}(H)} \, \Psi   = 0 .
\]
\end{quote}  
Then the next theorem follows.
\begin{quote} $\; $ \\
\textbf{Theorem I  (\cite{Arai00} ; Theorem 1.2)} \\Assume \textbf{(S.1)} and \textbf{(S.2)}. Then 
\[
\overline{  \{  E_0 ( H )  + \lambda \left| \frac{}{} \right.  \lambda \in  \sigma_{\ess} (S) \backslash \{ 0  \} }      \; \subset \sigma_{\ess} (H) ,
\]
where $\overline{J} $ denotes the closure of  $J \subset \mbf{R}$.
\end{quote}

$\;$ \\
$\;$ \\
{\large \textbf{Acknowledgments}} \\
It is  a pleasure to thank Professor  Tadayoshi Adachi  for giving  opportunities of talk in  seminar and  his  comments. This work is supported by JSPS grant 24$\cdot$1671. 

\end{document}